\documentclass[12pt,letterpaper]{article} 
\usepackage[margin=1in]{geometry} 

\usepackage{graphicx} 
\usepackage{amsmath}
\usepackage{amssymb}
\usepackage{natbib}
\usepackage{caption}
\usepackage{subcaption}
\usepackage{float}
\usepackage{soul}
\usepackage{hyperref}
\usepackage[table]{xcolor}
\usepackage[ruled]{algorithm2e}

\title{A Pseudo-Gradient Approach for Model-free Markov Chain Optimization}
\author{Nanne A. Dieleman\\ 
Department of Operations Analytics, Vrije Universiteit Amsterdam,\\
n.a.dieleman@vu.nl \\
Joost Berkhout\\
Department of Mathematics, Vrije Universiteit Amsterdam,\\
joost.berkhout@vu.nl
\\
Bernd Heidergott\\
Department of Operations Analytics, Vrije Universiteit Amsterdam,\\ 
b.f.heidergott@vu.nl \\
}
\date{ }

\begin{document}

\maketitle

\begin{abstract}

We develop a first-order (pseudo-)gradient approach for optimizing functions over the stationary distribution of discrete-time Markov chains (DTMC). We give insights into why solving this optimization problem is challenging and show how transformations can be used to circumvent the hard constraints inherent in the optimization problem. The optimization framework is model-free since no explicit model of the interdependence of the row elements of the Markov chain transition matrix is required. Upon the transformation we build an extension of Simultaneous Perturbation Stochastic Approximation (SPSA) algorithm, called stochastic matrix SPSA (SM-SPSA) to solve the optimization problem. The performance of the SM-SPSA gradient search is compared with a benchmark commercial solver. Numerical examples show that SM-SPSA scales better which makes it the preferred solution method for large problem instances. We also apply the algorithm to the maximization of web-page rankings in web-graphs based on a real-life data set. As we explain in the paper, when applying a first-order gradient search one typically encounters a phenomenon which we call ``infliction points," that is, jumps in the optimization trajectories between periods of almost stationary behavior that slow down the optimization. We propose a heuristic for avoiding such infliction points and  present a metastudy on a wide range of networks showing the positive effect of our heuristic on the convergence properties of SM-SPSA gradient search.\\

\noindent{\bf Keywords:} optimization, discrete-time Markov chain, pseudo-gradient method, stochastic approximation, stochastic matrix\\


\end{abstract}

\section{Introduction}

The analysis, design, and optimization of networks, such as power grid, communication, computer, queueing, social, and transportation networks, are of vital economic and social importance. Formally speaking, these networks can be modeled as directed weighted graphs, and the dynamic relations captured by these graph models can be analyzed via discrete-time Markov chains (DTMCs).

In this paper, we focus on optimizing a given general smooth (nonlinear) objective function $f$ over the stationary distribution of a Markov chain (representing a weighted graph) by adjusting the transition probabilities of the corresponding transition matrix. This problem is relevant in many practical settings, for example, for optimizing the page rank of a website and building or removing roads in a transportation network to reduce congestion; see, for example, \citet{de2008maximizing}, \citet{Salman2018}, \citet{Crisostomi2011}, and \citet{schlote2012traffic}. In this paper, we propose to solve this optimization problem using a first-order iterative optimization algorithm.

The considered problem has several challenging characteristics. First, feasible solutions are limited to stochastic matrices. This poses a challenge for first-order iterative methods since, in each step of the optimization algorithm, the candidate solution has to be a stochastic matrix. This constraint cannot be relaxed as the stationary distribution, and thus the objective function $f$, cannot be computed if the transition probabilities are changed so that the matrix is no longer stochastic. Second, the weights of the edges cannot be freely adjusted since, for example, increasing the transition probability from a node $i$ to a node $j$ requires decreasing the transition probabilities of (some) edges from $i$ to other nodes to ensure that the matrix remains stochastic. Third, while exact methods exist to solve the problem in the case of simple linear and quadratic objective functions $f$, the problem becomes more challenging by also considering general smooth (nonlinear) objective functions $f$. Moreover, exact methods do not scale well.

To overcome these challenges, we develop a transformation that maps transition matrices to a general matrix space. By doing so, the optimization problem can be treated as an unconstrained problem and becomes model-free, i.e., our approach does not require any parameterization of the transition probabilities. This has the additional benefit that the developed algorithm can be used for problems where there is no natural parameterization, which is, for instance, the case in social and web-graph networks. Moreover, if a parameterization model is used, the found solution may depend on this choice, whereas this is not the case with our model-free approach. For the actual optimization, we introduce a gradient descent algorithm elaborating on first-order (pseudo-)gradient methods. More specifically, we will introduce an adaptation of the simultaneous perturbation stochastic approximation (SPSA), see \citet{Spall2003}, for optimization in the context of Markov chains. As we show in the paper, the algorithm, called stochastic matrix SPSA (SM-SPSA), scales well and avoids calculating the gradient analytically, which becomes complex and computationally inefficient with increasing network size.

As we will illustrate in this paper, the gradient descent approach exhibits non-standard behavior, which we believe is of interest in its own right. In particular, the algorithm trajectories show ``infliction points'', which makes developing termination criteria challenging and slows down the optimization. We will overcome this problem by providing a heuristic for choosing the initial matrix, which leads to an efficient gradient-based optimization search. With this, we contribute to the development of stochastic approximation optimization methods.

In conclusion, this paper makes the following contributions:
\begin{itemize}
    \item A framework is introduced for performing model-free optimization of general smooth (nonlinear) objective functions over the stationary distribution of DTMCs. Our framework also allows for the optimization of functions of other DTMC measures (such as the Kemeny constant or deviation matrix), although the focus is on functions over the stationary distribution.
    \item We introduce stochastic matrix SPSA (SM-SPSA), which is an extension of SPSA for dealing with stochastic matrix constraints. 
    \item The phenomenon of infliction points is identified, and through numerical examples, it is shown that infliction points cannot be ignored when performing first-order optimization of stochastic matrices.
    \item A heuristic is introduced to reduce the impact of infliction points on the optimization.
\end{itemize}

The paper is structured as follows. Section~\ref{sec:problem_setting} formally introduces the problem, and a review of the relevant literature can be found in Section~\ref{sec:literature_review}. Section~\ref{sec:solver_approach} provides a mathematical formulation that can be solved with standard exact solvers for linear and quadratic objective functions $f$, and our first-order optimization approach is developed in Section~\ref{sec:pseudo_gradient}. The phenomenon of infliction points is addressed in Section~\ref{sec:infliction_points}: we introduce the concept of infliction points more properly, discuss how these influence the optimization, and present possible remedies. In Section~\ref{sec:numerical_experiments}, we benchmark the quality of the proposed SM-SPSA algorithm with that of a standard solver. We also discuss a real-life inspired web page ranking optimization problem that is too complex to solve with a standard solver and present the results obtained with the SM-SPSA algorithm. We conclude with a discussion, a conclusion, and ideas for further research in Section~\ref{sec:discussion_conclusion}.
\section{Problem Setting} \label{sec:problem_setting}


Let $ {\cal G} = ( S , E  , P)$ be a weighted directed graph, where $S=\{0,1,...,N-1\}$ is the set of nodes/states, $ E \subseteq S^2 $ the set of edges, and $ P_{m n } $ the weight of edge $ (m , n) $. Although we allow the edges $ ( m , n ) \in E$ to have zero weight, the weight of nonexistent edges $ ( m , n ) \not \in E$ is set to zero by definition. Throughout this paper, we assume that the weight matrix is stochastic, i.e., $ P_{ m n } \geq 0 $, for all $m , n \in S $, and $ \sum_{n \in S} P_{m n } =1 $ for all $ m \in S$; where we write $ A_{ m n } =(A)_{m n }$ for entry $ ( m , n ) $ of some matrix $ A $ unless it causes notational confusion. We apply the same notational convention to vectors. Consequently, $P$ is a transition probability matrix of a random walk on $(S, E)$, the location of which is described by a Markov chain with $P$. More formally, let $X_0 \in S$ denote the initial location of the random walk and $X_t \in S$ its location after $t$ transitions; then the discrete-time Markov chain $\{X_t, t\geq 0\}$, with $ \mathbb{P} ( X_{t+1} = n | X_t = m ) = P_{m n}$, for all $ m , n\in S $, models a random walk on $ ( S, E )$ with transition probability matrix $ P $. We call $ {\cal G} = (S, E, P)$, with $ P $ a Markov chain, a Markov graph. We will also use the terms stochastic (weight) matrix and Markov chain as synonyms. We assume

\begin{itemize}
    \item [{\bf (A1)}]
$ {\cal G} $ is strongly connected so that the Markov chain that represents the weight matrix is irreducible and aperiodic.
\end{itemize}
We denote by $\mathcal{P}$ the set of stochastic matrices that are irreducible and aperiodic Markov chains on $ ( S, E )$. For $ P \in \mathcal{P}$, the stationary distribution of $ P$ is uniquely defined and is denoted by $\pi ( P )$. Moreover, we let $C$ denote a binary matrix, where $ C_{ m n } = 1 $ indicates the edges $ ( m , n )\in E$ of which the weight can be adjusted in optimization, and $ C_{ m n } = 0 $ otherwise.  

For a given initial Markov chain $ P_0 \in \mathcal{P} $ and binary matrix $ C$, we denote by $\mathcal{P}( P_0 , C)$ the feasible subset of $\mathcal{P}$ that contains all stochastic matrices $Q$ with $ Q_{ m n } = (P_0)_{ m n } $ in case $ C_{ m n } = 0$. The initial Markov chain $P_0$ represents the status quo that will be optimized. Throughout the paper, we impose the following:
\begin{itemize}
    \item [{\bf (A2)}]
The initial Markov chain $ P_0 $ and the matrix $ C $ are fixed and given.
\end{itemize}
We will use Assumption {\bf (A2)} to simplify the notation by suppressing $ P_0 $ and $ C $ where this does not confuse, for example, we may write $\mathcal{P}$ instead of $\mathcal{P}( P_0 , C)$. 

Let $ f $ be an objective defined for probability vectors on $ S$. By {\bf (A1)}, $ \pi ( P ) $ is uniquely defined, and we consider the problem of maximizing (without loss of generality) $ f $ over $ \pi ( P ) $ for $ P \in \mathcal{P}(P_0 ,C)$. That is, simply denoting $ f ( \pi ( P ) ) $ by $ f_\pi ( P )$, we consider
\begin{equation} \label{optimisation_problem}
    \max_{P \in \mathcal{P}(P_0 ,C) } f_\pi (P)
\end{equation}
for $ P_0 \in \mathcal{P}$ and $C$ given (see {\bf (A2)}). Note that the $ P \in {\mathcal{P}(P_0 ,C) }$ for which $ \pi ( P )$ solves the above problem is typically not unique, as a stationary distribution does not uniquely determine a Markov chain.

Solving \eqref{optimisation_problem} numerically poses the challenge of encoding the feasible set $ \mathcal{P}(P_0 ,C)$ in an efficient way. For this, we impose the following condition:
\begin{itemize}
    \item [{\bf (A3)}] 
    For $ P \in \mathcal{P}( P_0 , C)$, the weights of the adjustable entries are bounded from below by some small constant $ \gamma > 0$. 
\end{itemize}
Assumption {\bf (A3)} ensures that $\pi(P)$ exists for all $ P \in \mathcal{P}( P_0 , C)$ (i.e., Markov chain $P$ remains irreducible and aperiodic to enforce Assumption {\bf (A1)}).

For ease of presentation we focus on functions of the stationary distribution because of the relevance of this setting in applications.
However, the developed framework applies to generic functions of Markov chains as well. Indeed, if the existence of the stationary distribution is of no concern, our framework can be straightforwardly adapted to the multichain case (with possible transient states), in which case condition {\bf (A3)} can be dropped. This allows one to optimize functions of the ergodic projector or the deviation matrix, for example, to optimize the structural ranking methodology from \citet{berkhout_ranking_2018}. Alternatively, one can focus on optimizing a function of the Kemeny constant, which can be seen as a measurement of chain connectivity \citep{Berkhout2016a}.

\section{Literature Review}\label{sec:literature_review}

The research presented in this paper touches on two major streams in the theory of Markov chain optimization: Markov graph optimization (Section~\ref{sec:literature_review_markov_chain_optimization}) and first-order methods (Section~\ref{sec:literature_review_first_order_pseudo_gradient_methods}).

\subsection{Markov Graph Optimization} \label{sec:literature_review_markov_chain_optimization}

A rich literature exists on the problem of adding or removing edges in Markov graphs to optimize objectives of the Markov chain; see, for example, \citet{Csáji2010}, \citet{Fercoq2013}, \citet{Olsen2014}, \citet{Rosenfeld2016}, and \citet{Salman2018}.
In this paper, we consider a generalization of this problem by allowing for continuously adjustable edge weights. Moreover, we allow for controlling the set of edges whose weights can be changed.

Allowing for continuously adjustable edge weights links our approach to the rich literature on perturbation analysis of Markov chains. In this stream of research, a parameterization of the entries of the Markov chains is assumed, and then the optimization is carried out with respect to the model parameter; see, for example, \citet{cao2008stochastic}. For work related to queuing theory, parameterized models are typically available where the considered parameters are arrival rates or service rates. For these types of problems, closed-form solutions for the gradient of the stationary distribution can be found in the literature; see, for example, \citet{heidergott2003taylor} and \citet{meyer1994sensitivity}. However, in the analysis of (general) networks, such as social networks or web graph applications, there is typically no natural parameterization available. To overcome this problem, one typically postulates {\em ex-ante} a functional relationship between the row elements of the weight matrix of the network; see \citet{caswell2019sensitivity} and \citet{Berkhout2016a}. In contrast, we take a different approach by considering transformations, which means that it is unnecessary to assume a functional relationship between the row elements, making our approach model-free.

\subsection{First-Order and Pseudo-Gradient Methods} \label{sec:literature_review_first_order_pseudo_gradient_methods}

First-order and pseudo-gradient methods are versatile optimization tools that can be used in a variety of settings. 
Pseudo-gradient methods use an estimate/proxy for the gradient. Examples are Simultaneous Perturbation Stochastic Approximation (SPSA) \citep{Spall2003} and Smoothed Functional (SF) algorithms such as the Gaussian Based SF Algorithm (GSFA) \citep{Bhatnagar2013}.

Applying gradient methods to optimize the transition probabilities of a Markov chain leads to a constrained optimization problem, and projection is needed to ensure that the optimization parameters remain feasible. See, for example, \cite{bertsekas1976goldstein} for an explanation of projected gradient search. Although a projection on the probability simplex is available \citep{Wang2013}, applying this projection row-wise in combination with the analytical gradient to optimize over $ \mathcal{P}$ renders the gradient descent numerically inefficient. Mirror descent, introduced by \citet{Nemirovsky1983}, offers an alternative to projection for constrained optimization,
and we refer to \cite{d2021stochastic} for a mirror descent algorithm to find the stationary distribution of a Markov chain. Unfortunately, mirror descent requires a nonlinear transformation of the gradient, which rules out the possibility of combining it with randomized pseudo-gradient methods like SPSA. However, using these pseudo-gradient methods has considerable advantages, such as computational efficiency, over analytically computing the gradient, which is especially important for larger networks. To avoid projection or any nonlinear transformation of the gradient, we introduce a transformation of the parameter space $\mathcal{P}$ and apply the (pseudo-)gradient search in the transformed parameter space. 


\section{An Exact Approach} \label{sec:solver_approach}


In this section, we formulate a mathematical formulation for optimization problem \eqref{optimisation_problem}. This mathematical formulation can be solved by a solver for specific forms of the objective function $f$. 
With $ {\bf (A1)} $ to $ {\bf (A3)} $ in force, the mathematical formulation is given by:
\begin{align}
    \text{max}\quad & f_\pi( P ) & \label{Mathematical_formulation} \\
    \text{s.t.}\quad &\sum_{n \in S} P_{mn} = 1&& \forall m \in S \label{Mathematical_model_2}\\
                &  (\pi(P))_m = \sum_{n \in S} (\pi(P))_n  \, P_{nm}&& \forall m \in S \label{Mathematical_model_3}\\
                & \sum_{m \in S} (\pi (P))_{m} = 1&& \label{Mathematical_model_4}\\
                & P_{mn} = \bigl(P_0\bigr)_{mn} && \forall m,n \in S: C_{mn} = 0 \label{Mathematical_model_5}\\
                & (\pi (P))_m \geq 0&& \forall m \in S \label{Mathematical_model_6}\\
                & P_{mn} \in [ \gamma, 1] &&\forall m,n \in S \label{Mathematical_model_7}
\end{align}

Constraint \eqref{Mathematical_model_2} ensures that all rows of the transition matrix sum up to 1. Constraint \eqref{Mathematical_model_3} constitutes the balance equations. Constraint \eqref{Mathematical_model_4} ensures that the stationary distribution sums up to 1. Transition probabilities $P_{mn}$ that are not allowed to be adjusted will keep initial transition probabilities $\bigl(P_0\bigr)_{mn}$, which is specified by constraint \eqref{Mathematical_model_5}. Constraint \eqref{Mathematical_model_6} forces all stationary probabilities to be nonnegative. Together with constraint \eqref{Mathematical_model_4}, this ensures that they are probabilities. Lastly, \eqref{Mathematical_model_7} enforces ({\bf A3}). 

This mathematical model can be solved with state-of-the-art commercial software, such as Gurobi \citep{Gurobi2023}, if the objective function $f$ is linear or quadratic. 
However, using Gurobi or similar solvers has multiple drawbacks. First, this solution approach does not scale well. In the worst case, its time complexity is exponential since it is a non-convex continuous model. We will illustrate this with numerical experiments in Section~\ref{sec:numerical_experiments}. Second, no general smooth (nonlinear) objective functions can be straightforwardly solved with this method, while these types of objective functions often arise in practice. Lastly, the constraints of the mathematical model are satisfied up to a user-defined tolerance, which means that, in some cases, the final solution may not meet all constraints due to this tolerance. Moreover, we observed that in some experiments, this tolerance was not satisfied, probably due to numerical issues of the solver used (Gurobi version 10).

In the subsequent section, we will present a pseudo-gradient approach that can also solve the provided mathematical model but does not have the three drawbacks discussed previously.

\section{A Pseudo-Gradient Approach} \label{sec:pseudo_gradient}

In this section, we introduce the framework that can solve optimization problem \eqref{optimisation_problem} by using an iterative pseudo-gradient approach based on stochastic approximation. It consists of a combination of an iterative first-order method and transformations. We first introduce the general stochastic matrix optimization framework for a first-order method in Section~\ref{sec:first_order_method}. Subsequently, we show how SPSA can be used in combination with the framework in Section~\ref{SPSA_formulation}.

\subsection{The General Stochastic Matrix Optimization Framework}\label{sec:first_order_method}

First-order methods are gradient search-based optimization methods. In general terms, let $ {\cal X } $ denote the feasible set, $ f $ some sufficiently smooth differentiable function defined on $ {\cal X}$, then a first-order method to solve $ \max_{\vartheta \in {\cal X}  }  f ( \vartheta ) $ takes the form
\begin{equation}\label{Standard_SA}
    {\vartheta}^{(i+1)} = \Pi_{ {\cal X}} \Big ( {\vartheta}^{(i)}+\epsilon
     {G}({\vartheta}^{(i)}) \Big ), 
    \quad i \geq 0 , 
\end{equation}
where $i$ represents the iteration number, $ G (\vartheta) $ the gradient of $ f $ at $\vartheta$, $ \Pi_{ {\cal X}} $ the projection on $ {\cal X} $, $\epsilon$ the constant step size, and $ \vartheta^{(0)}$ the initial value of the algorithm. Under appropriate conditions, the limit point of the sequence $\{ \vartheta^{(i)} \}$ is a stationary point of $ f $, so that under appropriate smoothness of $ f $ this point is the location of a (local) maximum of $ f$ over $ \cal X$, see \cite{ Bert} for details. The algorithm in \eqref{Standard_SA} extends to the setting where $ G $ is a pseudo-gradient, i.e., when $ G $ is an approximation of the true gradient \citep{Spall2003}.

There is no straightforward way to apply the first-order method in \eqref{Standard_SA} to the optimization problem in \eqref{optimisation_problem}. To see this, let $ \vartheta $ denote the vector containing the adjustable elements of the Markov chains; that is, each $ \vartheta$ represents all entries $ P_{m n } $ with $ C_{m n } =1 $, and denote the Markov chain under $ \vartheta $ by $ P ( \vartheta )$. First, $ \partial P ( \vartheta ) / \partial \vartheta$ is not well-defined as perturbing $ \vartheta$ requires assuming a model of how to redistribute the change in $ \vartheta$ to keep the perturbed matrix stochastic, see \citet{caswell2019sensitivity} and \citet{Berkhout2016a}, which is contrary to our aim of developing a model-free approach. 
Second, to evaluate the stationary distribution and thus the objective function, $ P( \vartheta ) $ must be a stochastic matrix, and $ P ( \vartheta ) \in {\mathcal{P}}( P_0 , C )$ is thus a hard constraint that cannot be softened during optimization. This rules out, for example, constraint-violation penalty methods. Theoretically, a projection $\Pi_{\cal X}$ on the set $ \cal X $ of admissible parameters could be used. However, computing the gradient analytically is computationally inefficient  and therefore we approximate the gradient by a perturbation method, such as SPSA or GSFA (see \citet{Spall2003} and \citet{Bhatnagar2013}). Unfortunately, these perturbation methods automatically violate the stochastic matrix constraint due to the use of random perturbations when approximating the gradient.

To overcome this drawback, we propose to transform the optimization problem in \eqref{optimisation_problem} into an unconstrained problem on $\mathbb{R}^{N \times N}$. To do so, we transform the matrix $P_0 \in \mathcal{P}$ to a matrix $\Theta^{(0)} \in \mathbb{R}^{N \times N}$ with a function to be introduced later in the text, and let the first-order algorithm in \eqref{Standard_SA} run in $\mathbb{R}^{N \times N}$, where we update a candidate matrix $ \Theta^{(i)} $ in iteration $i$ element-wise via \eqref{Standard_SA} for the adjustable entries only. In order to apply \eqref{Standard_SA}, we need to transform matrices in the optimization space $\mathbb{R}^{N \times N}$ back to stochastic matrices in the space $\mathcal{P} (P_0, C)$. We do so by the transformation $T(\Theta)$, to be introduced later. A key benefit of the transformation is that it enables the use of perturbation methods to approximate the gradient since the gradient is computed in the \textit{unconstrained} optimization space $\mathbb{R}^{N \times N}$. Note that the outer projection in \eqref{Standard_SA} can also be ignored for the same reason. The transformation does not introduce new stationary points, which is an important characteristic when it is applied in combination with an optimization algorithm, which is the case in this paper.

A matrix $\Theta \in \mathbb{R}^{N \times N}$ can be transformed back to a stochastic matrix in the space ${\mathcal{P}}( P_0 , C )$ in three steps. The first two steps ensure that a stochastic matrix is obtained. The third step is a row-scaling operation that takes into account the binary adjustment matrix $C$. More specifically, the transformation takes the following steps:

\begin{itemize}
\item \textbf{Step 1:} The adjustable entries where $C_{mn}=1$ of matrix $\Theta$ are (entry-wise) transformed by a smooth function $\Phi^{\rm{entry}}(\cdot)$: $\mathbb{R} \rightarrow \mathbb{R}^{+}$, defined as the standard logistic function:
\begin{equation}
    (\Phi^{\rm{entry}}(\Theta))_{mn} = 
    \frac{1}{1+e^{-\Theta_{mn}}}  \qquad\qquad\qquad \forall m,n \in S: C_{mn} = 1.
\end{equation}
We denote the matrix output of this function by $V$. 
\item \textbf{Step 2:} All rows are normalized so that the sums of the rows of the adjustable entries are equal to 1 by the function $\Phi^{\rm{row}}(\cdot)$, given by
\begin{equation}
    (\Phi^{\rm{row}} (V))_{mn} = \frac{V_{mn}}{\sum_{k=1}^N 1_{\{C_{mk}=1\}} V_{mk}} \qquad\qquad\qquad \forall m,n \in S: C_{mn} = 1.
\end{equation}
We denote the matrix output of this function by $U$. 
\item \textbf{Step 3:} The rows are scaled so that all the rows sum up to 1, adjusting for the entries that cannot be adjusted according to matrix $C$, by the function $\Phi^{\rm{scale}}(\cdot)$, defined by
\begin{equation} \label{function_scaling_operation}
    (\Phi^{\rm{scale}}(U))_{mn} = 
    U_{mn} \left(1-\sum_{k=1}^N 1_{\{C_{mk}=0\}} U_{mk}\right) \qquad\forall m,n \in S: C_{mn} = 1.
\end{equation}
\end{itemize}
The complete transformation $T$ is then given by
\begin{equation}
    T(\Theta) = \Phi^{\rm{scale}}\circ 
 \Phi^{\rm{row}} \circ \Phi^{\rm{entry}} ( \Theta ) 
 .   
\end{equation}
Other functions than the standard logistic function are also valid for $\Phi^{\rm{entry}}$, and extensive numerical studies showed that these functions should be bounded and symmetric because of numerical stability.

\begin{figure}[H]
    \centering
    \centerline{
    \includegraphics[scale=0.55]{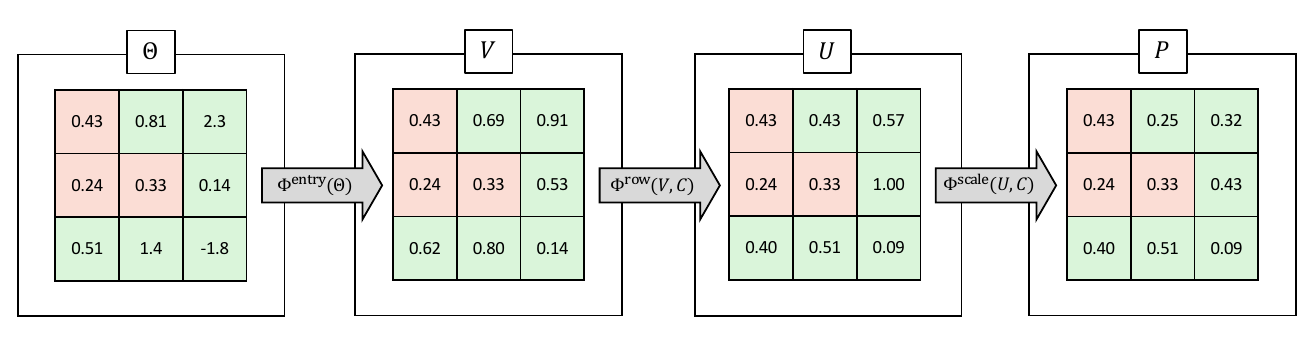}
    }
    \caption{Example of the complete transformation $T(\Theta)$ that transforms matrix $\Theta$ into stochastic matrix $P \in \mathcal{P}$. The red cells are non-adjustable (its $C$-entries are 0), and the green cells are adjustable (its $C$-entries are 1).}
    \label{fig:Transformation_calculation}
\end{figure}

Figure~\ref{fig:Transformation_calculation} illustrates the complete transformation for a small $3 \times 3$ randomly generated matrix $\Theta \in \mathbb{R}^{N \times N}$. The green and red boxes represent the entries that can and cannot be adjusted according to the matrix $C$, respectively.

For initializing the algorithm, the initial transition matrix $P_0$ is transformed into a matrix $\Theta^{(0)}$ in $\mathbb{R}^{N \times N}$ by a function $(\Phi^{\rm{entry}})^{-1}$. This function $(\Phi^{\rm{entry}})^{-1}$ is simply the inverse of the function $\Phi^{\rm{entry}}$, and in the case of a standard logistic function equal to
\begin{equation}
    \bigl((\Phi^{\rm{entry}})^{-1}(P)\bigr)_{mn} = \ln\left(\frac{P_{mn}}{1-P_{mn}}\right) \qquad \forall m,n \in S: C_{mn} = 1.
\end{equation}
It should be noted that a problem arises when the function $ (\Phi^{\rm{entry}})^{-1}$ is applied to the entries where $P_{mn} =0$ or $P_{mn}=1$, for all $m, n \in S$. This problem can be solved by changing zero values to $\gamma$ and values of one to $1-\gamma$ when applying $ (\Phi^{\rm{entry}})^{-1}$ with $\gamma > 0$ small. For example, one could choose $\gamma$ to be slightly smaller than the smallest nonzero entry. However, it should be noted that the use of small $\gamma$ can lead to infliction points, a phenomenon that is discussed further in Section~\ref{sec:infliction_points}. Generally, it is advisable to keep numerical precision in mind due to floating-point arithmetic. It is also interesting to note that with the choice for this particular function for $\Phi^{\rm{entry}}$, explicitly enforcing Assumption~{\bf (A3)} is not necessary. This is because the values 0 and 1 are only asymptotically reachable, which means that the Markov chain will remain irreducible during optimization. It might be the case that the actual optimal solution has transition probabilities of value 0 and/or 1. To find this solution, one must manually check whether the transition probabilities close to 0 or 1 should exactly be 0 or 1. This is also the case for the mathematical model.

\begin{figure}[H]
    \centering
    \centerline{\includegraphics[scale=0.62]{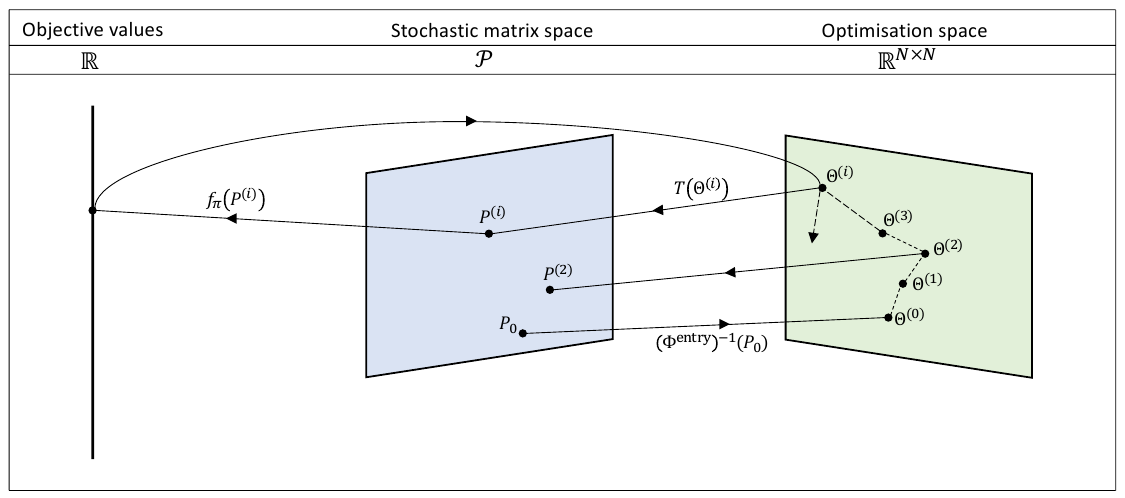}}
    \caption{Illustration of the proposed framework. We start the optimization with transition matrix $P_0$ (=$P^{(0)}$), which is transformed to the matrix ${\Theta}^{(0)}$. We then use the first-order method in \eqref{Standard_SA} (without projection) to calculate ${\Theta}^{(1)}$, which corresponds to $P^{(1)}$ and objective value $f_\pi (P^{(1)})$; and so forth.} 
    \label{fig:Transformation_illustration}
\end{figure}

To conclude, the framework is as follows. First, we transform the initial matrix $P_0$ by the function $(\Phi^{\rm{entry}})^{-1}$, which produces the matrix $ \Theta^{(0)} \in \mathbb{R}^{N \times N}$. We then apply the first-order method in \eqref{Standard_SA} (without projection), where each $ \Theta^{(i)}$ can be transformed back to a Markov chain $ P^{(i)}$ by transformation $T$, so that $ f_\pi $ can be evaluated and the gradient approximated. See Figure~\ref{fig:Transformation_illustration} for an illustration of the framework. In the next section, we show how the gradient can be approximated by combining this framework and SPSA.

\subsection{The Stochastic Matrix SPSA Algorithm} \label{SPSA_formulation}


Simultaneous Perturbation Stochastic Approximation (SPSA) was developed by \citet{Spall2003}. It is a popular stochastic optimization algorithm, as it scales well and is easy to implement. In this section, we discuss an adaptation of the SPSA algorithm by using the framework introduced in Section~\ref{sec:first_order_method} to solve problem \eqref{optimisation_problem}. We call the algorithm ``Stochastic Matrix Simultaneous Perturbation Stochastic Approximation" (SM-SPSA).

The algorithm's input is an initial transition matrix $P_0$ and a binary adjustment matrix $C$. The first step is to apply the inverse transformation $(\Phi^{\rm{entry}})^{-1}$ to the matrix $P_0$ to obtain a matrix ${\Theta}^{(0)} \in \mathbb{R}^{N \times N}$. Next, the iterative component of the algorithm starts. The fixed gain size recursive stochastic approximation formula in \eqref{Standard_SA} (without projection) is used. The gradient estimate $G^{(i)}$ at iteration $i$ is calculated by
\begin{equation} \label{Matrix_SPSA_Gradient}
    \Bigl({G}^{(i)}\Bigr)_{mn} =  C_{mn} \frac{f^{+, (i)} - f^{-, (i)}
    }{2\eta_i(\Delta^{(i)})_{mn}},
\end{equation}
with
\[
f^{\pm, (i)} = f_\pi \left(T\left({\Theta}^{(i)} \pm  \eta_i\left(\Delta^{(i)} \odot C\right)\right)\right),
\]
and where: 
\begin{itemize}

\item $\Delta^{(i)}$ is a random perturbation matrix with i.i.d.
entries $ \mathbb{P} ( (\Delta^{(i)} )_{ m n } =- 1 ) = 0.5 = 
\mathbb{P} ( (\Delta^{(i)} )_{ m n } = 1 )  $,

\item $ \Delta^{(i)} \odot C $ denotes the Hadamard product (i.e., element-wise product) of $ \Delta^{(i)} $ and $ C $,

\item $\eta_i = \frac{1}{i+1}$ is a positive and decreasing sequence scaling the finite perturbations.

\end{itemize}
Provided that $ f_\pi ( T ( \Theta ) ) $ is a smooth mapping of $ \Theta   $ with bounded third-order derivatives, the algorithm is known to converge to a stationary point of $ f_\pi ( T ( \Theta ) ) $ for $ \epsilon$ sufficiently small. The convergence theory was originally developed for decreasing gain size; see \cite{spall1992multivariate}, and for the limit theory for fixed gain size we refer to \cite{gerencser2000spsa,gerencser2001mathematics}.

Algorithm \ref{Matrix_SPSA_algorithm} presents the algorithm step by step. The SM-SPSA algorithm code is available at \url{https://github.com/nanned/SM-SPSA} and a website with set-up guidelines, documentation and explanations at \url{https://nanned.github.io/SM-SPSA/}. For $I$ iterations, the algorithm has worst-case complexity $O(IN^3)$, attributed to the $O(N^3)$ worst-case complexity of calculating the stationary distribution.

\begin{algorithm}[H]
\caption{The SM-SPSA algorithm to solve Problem~\eqref{optimisation_problem}} \label{Matrix_SPSA_algorithm}
\KwData{Initial transition matrix $P_0$, binary adjustment matrix $C$, $\epsilon$, total number of iterations $I$.}
\KwResult{Optimised transition matrix}
${\Theta}^{(0)}$ = $(\Phi^{\rm{entry}})^{-1}(P_0)$\;
\For{i = 0 \KwTo I - 1}{
Generate perturbation matrix $B^{(i)} = \eta_i(\Delta^{(i)} \odot C)$ \;
$\Theta^{+, (i)} = {\Theta}^{(i)} + B^{(i)} $ \;
$\Theta^{-, (i)} = {\Theta}^{(i)} - B^{(i)} $ \;
Calculate the stationary distributions of $\Theta^{+, (i)}$ and $\Theta^{-, (i)}$ by using the transformation $T$: $\pi^{+, (i)} = \pi(T(\Theta^{+, (i)}))$ and $\pi^{-, (i)} = \pi(T(\Theta^{-, (i)}))$\;
Calculate the gradient proxy ${G}^{(i)}= \frac{f(\pi^+_{i})-f(\pi^-_{i})}{2\eta_i\Delta_i} \odot C$\;
Update ${\Theta}^{(i+1)} = {\Theta}^{(i)}+\epsilon {G}^{(i)} $\;
}
\end{algorithm}

\section{The Price of Transformations: Infliction Points} \label{sec:infliction_points}

Using a transformation to transform a constrained optimization problem into an unconstrained optimization problem will come at a price. In our case, the values 0 and 1 are only asymptotically reachable, for example. Moreover, we encountered an interesting phenomenon, which we call infliction points. The phenomenon will be discussed in detail in Section~\ref{sec:infliction_points_phenomenon_of_infliction_points}, and in Section~\ref{sec:infliction_points_remedies}
we introduce heuristics for avoiding such infliction points.

\subsection{The Phenomenon of Infliction Points} \label{sec:infliction_points_phenomenon_of_infliction_points}
To explain these infliction points, we consider the smallest network where such an infliction point can arise, namely a network with 3 nodes. Let
\begin{equation}\label{eq:example_infliction_points}
    P_0 =
    \begin{pmatrix}
    0.001 & 0.001 & 0.998 \\
    0.998 & 0.001 & 0.001\\
    0.001 & 0.998 & 0.001
    \end{pmatrix},
    C = \begin{pmatrix}
    0 & 1 & 1 \\
    1 & 0 & 1\\
    1 & 1 & 0
    \end{pmatrix}
\end{equation}
be the initial transition matrix and the binary adjustment matrix, respectively. The objective is to maximize the stationary distribution of node 0. Three graphs that illustrate the results of running SM-SPSA on this network are shown in Figure~\ref{fig:infliction_points_example}. The used SM-SPSA parameters are provided in Appendix~\ref{Appendix_test}. The upper graph shows the objective values for each iteration $i$. The middle graph shows the values of the adjustable transition probabilities $(P^{(i)})_{mn}$ $\forall i$ in the stochastic matrix space. The lower graph shows the run in the optimization space for each adjustable entry $(\Theta^{(i)})_{mn}$ $\forall i$.


Infliction points are the points where 
jumps in the optimization path occur after a relatively long period where no changes seem to occur. In Figure~\ref{fig:infliction_points_example}, this occurs at points (a) and (b). Preceding both points, it seems as if the optimization has converged to a specific solution, whereas after the infliction points, the optimization again seems to converge, but then to different solutions. Thus, the optimization switches abruptly between solutions at infliction points. This phenomenon is important, as it is challenging to decide when to stop the optimization based on the objective and the $(P^{(i)})_{mn}$ values.
This is best observed in the graph of the objective function value: it seems as if the optimization has converged in the flat segment between (a) and (b). However, at point (b), the objective function value suddenly jumps from $\pm$ 0.45 to $\pm$ 0.50. Standard stopping criteria based on the objective values and/or the $(P^{(i)})_{mn}$ values would therefore already stop the optimization before point (a) or in the constant segment between (a) and (b), whereas, in fact, the optimization has not converged yet. It is also interesting to observe that these sudden changes are not visible in the graph of the $(\Theta^{(i)})_{mn}$ values. Another, and perhaps most important, disadvantage of the presence of infliction points is that the optimization itself also requires more iterations before it converges because of these long periods of quasi-stationarity. This makes the overall optimization slower. 

\begin{figure}[H]
    \centering
    \centerline{
    \includegraphics[scale=0.45]{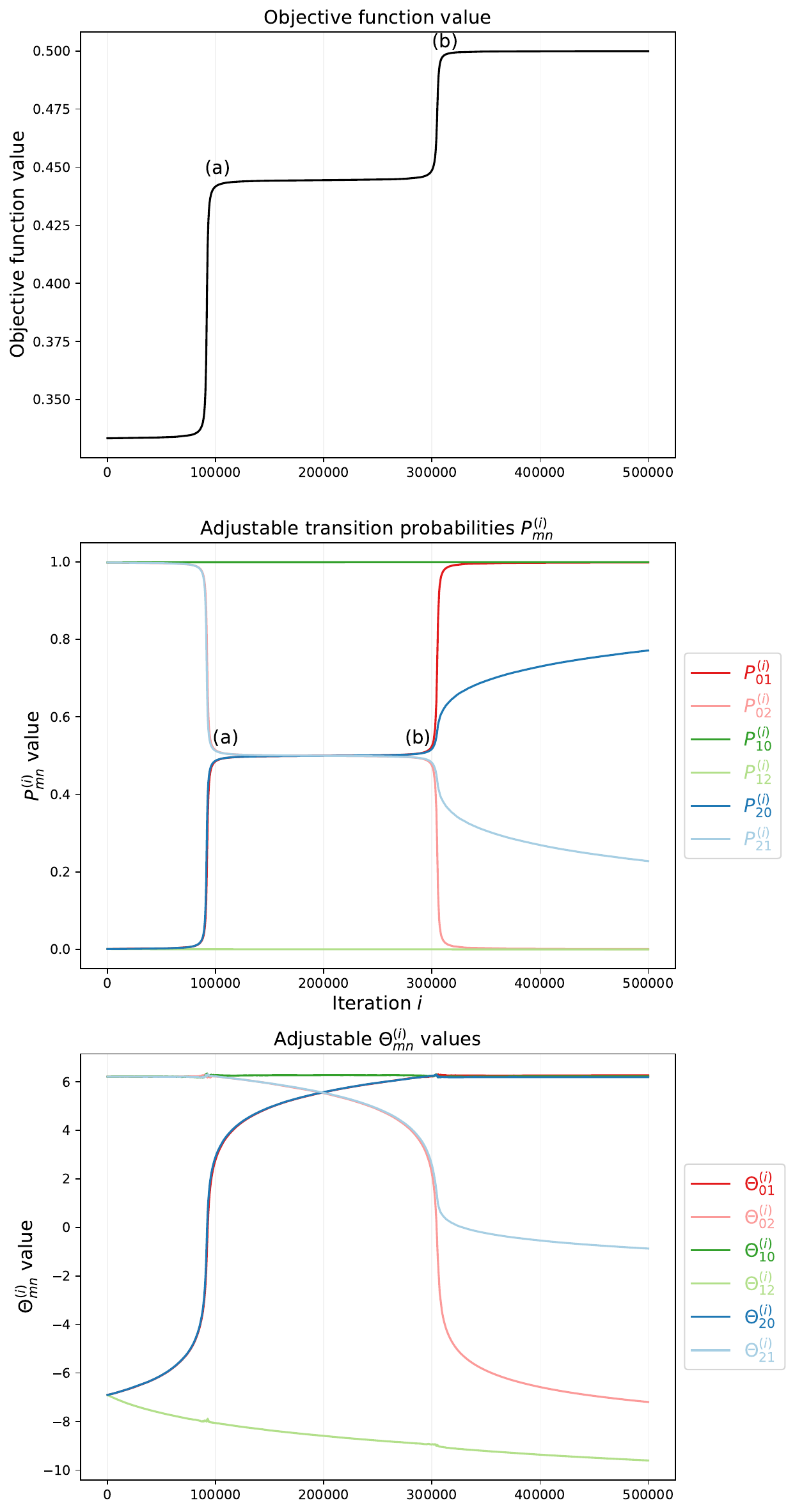}
    }
     \caption{Results of an optimization run for the network in \eqref{eq:example_infliction_points} where infliction points occur at points (a) and (b).}
    \label{fig:infliction_points_example}
\end{figure}

These infliction points arise due to small transition probabilities $(P_0)_{mn}$ combined with the transformation $T$. The transformation $T$ has a ``squashing" effect when applied to matrices in the optimization space $\mathbb{R}^{N \times N}$. This is inherent to the transformation, as in the optimization space, the entries can attain any real value, although in the stochastic matrix space $\mathcal{P}$, the entries are restricted to be nonnegative and have row sums of 1. This squashing effect makes it more difficult for a small $P_{mn}$ entry to grow to a larger value since a relatively much larger change in the entry $\Theta_{mn}$ is required. This can, for example, be observed with entry $P_{20}$ (corresponding to $\Theta_{20}$) in Figure~\ref{fig:infliction_points_example}. It should be noted that these infliction points are relevant in practice, as often real-world networks contain weak links, i.e., edges with small weights.

\subsection{Remedies for Infliction Points}\label{sec:infliction_points_remedies}
A first idea to deal with these infliction points may be to base the stopping criterion on the $({\Theta}^{(i)})_{mn}$ values in the optimization space $\mathbb{R}^{N \times N}$ rather than on the $({P}^{(i)})_{mn}$ values or the objective values. A common stopping criterion is then to stop the algorithm at iteration $i$ ($\geq R$) if $||{\Theta}^{(i)}-{\Theta}^{(i-j)}|| \leq \omega$ holds for all $1 \leq j\leq R$, for some small $\omega > 0$ and $R$ the number of iterations for which the algorithm must be stable \citep{Spall2003}. Although this stopping criterion ensures that the optimization is not prematurely stopped, it does not deal with the second disadvantage of infliction points, namely the stagnant optimization itself.

To overcome the problem of stagnant optimization, we propose a starting heuristic that deals with small values of $P_0$. Given an initial matrix $P_0$, we rearrange the mass of the adjustable entries of $P_0$ per row by evenly assigning the available mass to these entries. All adjustable entries then start row-wise at the same point. We then perform the optimization on this new starting matrix. We call this heuristic the ``centered mass" heuristic. Our hypothesis is that the centered mass heuristic results in a better objective value after $t$ iterations compared to the situation without this heuristic since no or at least fewer infliction points will occur. This indicates that the optimization is indeed faster with the centered mass heuristic. Appendix~\ref{Appendix_test} contains the graphs for running SM-SPSA on the network in \eqref{eq:example_infliction_points} with the centered mass heuristic; these graphs illustrate the profound positive influence that the centered mass heuristic can have on the optimization. 

We performed several statistical experiments with three different network sizes to verify our hypothesis. We test this hypothesis in a realistic setting where there is a limited computational budget. Let $p_t$ be the probability that the centered mass heuristic results in a higher objective value after $t$ iterations. The hypothesis then becomes:

\begin{align*}
    H_0(t): p_t\leq\frac{1}{2}\\
    H_1(t): p_t>\frac{1}{2}
\end{align*}

We generated instances of three different problem classes, namely ``small" networks (10 nodes, 500 instances), ``medium" networks (50 nodes, 10 instances), and ``large" networks (100 nodes, 10 instances) according to the sampling method described in Appendix~\ref{Appendix_sampling_networks}. The randomly generated binary adjustment matrices have zeroes on their diagonals to avoid changing self-loops. Otherwise, the optimization is trivial. For each network, the objective is to maximize the stationary distribution of a randomly selected node. We ran SM-SPSA on all instances of each problem class and counted the number of times the centered mass heuristic resulted in a higher objective value after $t \in \{1000, 2000, ...., I\}$ iterations. We then performed a population proportion test for each $t$ and each problem class. The parameters used to run the SM-SPSA algorithm can be found in Table~\ref{tab:statistical_experiment_parameters_SM-SPSA} in Appendix~\ref{Appendix_test}.

\begin{figure}[H]
    \centering
    \centerline{
    \includegraphics[scale=0.5]{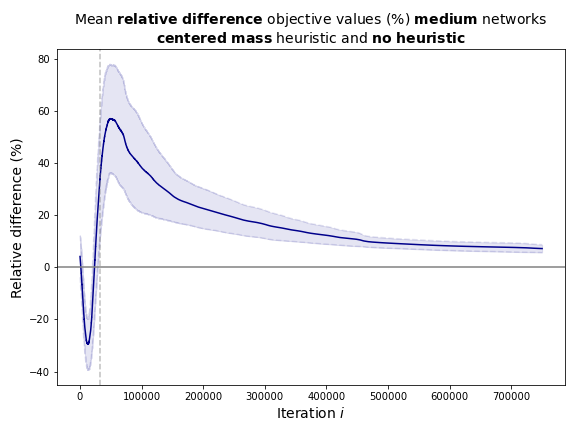}
    }
     \caption{Results for the medium-class networks.}
     \label{fig:mean_rel_diff_heuristic_medium}
\end{figure}

\begin{figure}[H]
    \centering
    \centerline{
    \includegraphics[scale=0.5]{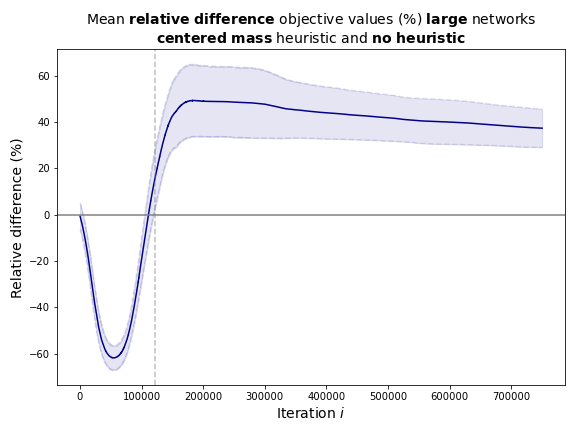}
    }
     \caption{Results for the large-class networks.}
     \label{fig:mean_rel_diff_heuristic_large}
\end{figure}

Figures~\ref{fig:mean_rel_diff_heuristic_small}, \ref{fig:mean_rel_diff_heuristic_medium} and \ref{fig:mean_rel_diff_heuristic_large} present the mean relative difference (\%) with 95\% confidence intervals between the objective values with and without using the centered mass heuristic during the optimization for the small-, medium- and large-class networks, respectively. For example, a mean relative difference of 10\% means that the found objective value by SM-SPSA with the centered mass heuristic is, on average, 10\% better than without the heuristic. The vertical line in each graph represents the iteration $t \in \{1000, 2000, ...., I\}$ from which iteration onward the null hypothesis is rejected. 

\begin{figure}[H]
    \centering
    \centerline{
    \includegraphics[scale=0.5]{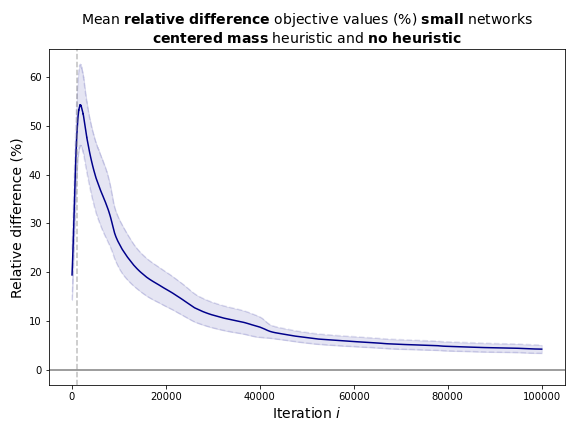}
    }
     \caption{Results for the small-class networks.}
    \label{fig:mean_rel_diff_heuristic_small}
\end{figure}

\begin{figure}[H]
    \centering
    \centerline{
    \includegraphics[scale=0.5]{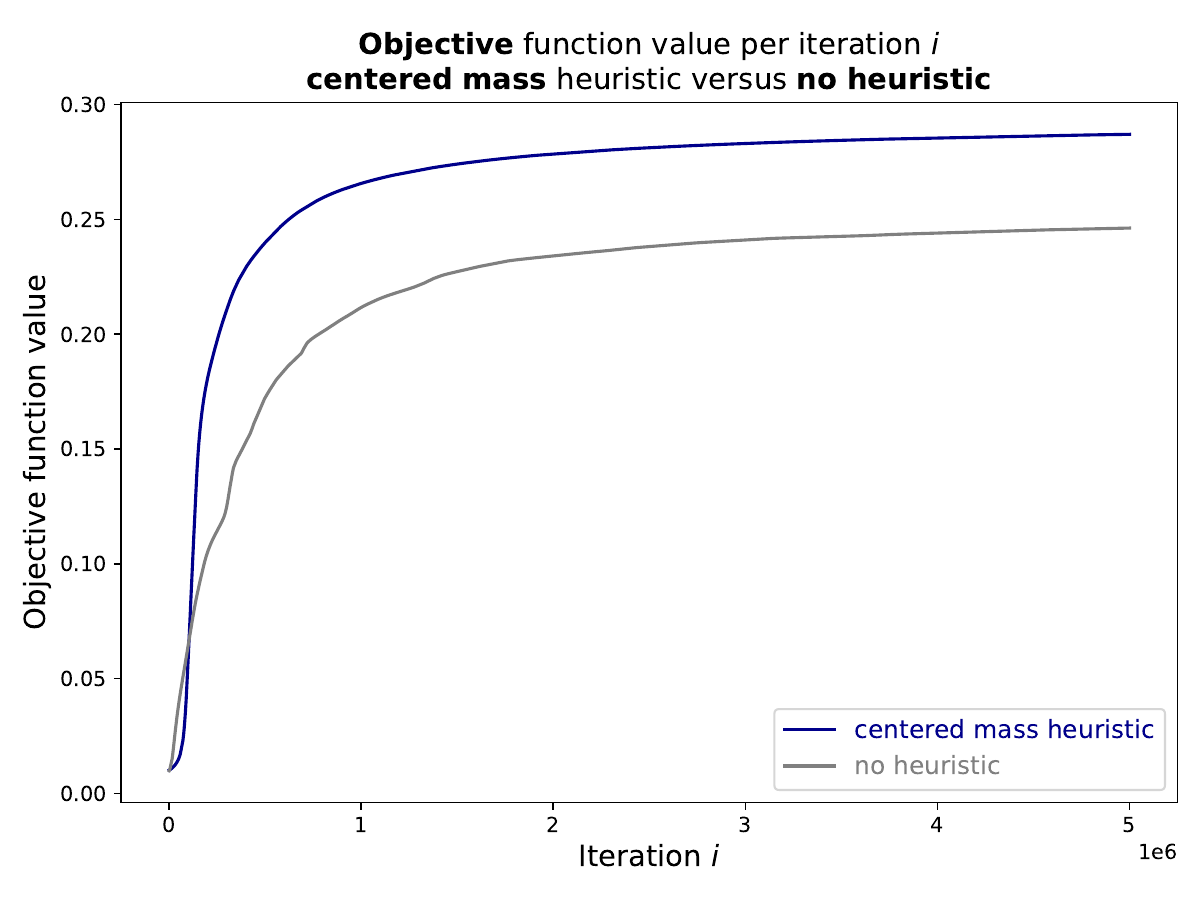}
    }
     \caption{The objective values during an optimization of a large-class network with and without using the centered mass heuristic.}
     \label{fig:long_run_centred_mass_heuristic}
\end{figure}

In the case of the medium- and large-class networks, we can observe a ``warm-up" period during which not using the heuristic is better than using the centered mass heuristic. The centered mass heuristic can be seen as an ``uninformed prior" that gives each adjustable transition probability an equal opportunity to increase or decrease in value during the optimization. On the other hand, not using the heuristic already ``favors" a certain (local) optimum, since some transition probabilities already have smaller or larger values than others. The optimization continues in the direction of this (local) optimum and may experience an almost standstill in the run-up to an infliction point, which slows down the optimization. The objective values are therefore higher with the centered mass heuristic after the warm-up period, resulting in rejected null hypotheses. The graphs show that these warm-up periods become longer as the networks grow. However, the positive influence of the centered mass heuristic also lasts longer for larger networks. To illustrate this, Figure~\ref{fig:long_run_centred_mass_heuristic} shows the objective values with and without the heuristic during optimization of a large-class network over 5.000.000 iterations using $\epsilon=0.1$. Using the heuristic indeed results in higher objective values after the warm-up period. The ``bumps" in the graph of the optimization without heuristic correspond to infliction points.

To summarize, as long as the computational budget allows to run the algorithm longer than a relatively small warm-up period, these experiments and tests show that the centered mass heuristic indeed results in significantly better objective values and thus a significantly faster optimization compared to not using the heuristic. In short, we recommend using the centered mass heuristic, except if only a limited number of iterations are used. Unless otherwise stated, we use the centered mass heuristic in the numerical experiments of Section~\ref{sec:numerical_experiments}.


We conclude this section by noting that the phenomenon of infliction points is not specifically a characteristic of the SM-SPSA algorithm but rather a characteristic of any first-order method based on the framework proposed in Section~\ref{sec:first_order_method}. For example, similar behavior would be observed if a finite-difference gradient estimation is used.

\section{Numerical Experiments} \label{sec:numerical_experiments}

In this section, we perform multiple numerical experiments to evaluate the performance and versatility of the SM-SPSA algorithm. In Section~\ref{sec:numerical_experiments_SM-SPSA_vs_Gurobi}, we will compare the solution quality of the SM-SPSA algorithm to that of the Gurobi solver for various relatively smaller and simpler problems. We will also compare the performance of both algorithms for a large network of 500 nodes. In Section~\ref{sec:numerical_experiments_hyperlink_optimization}, we will show the applicability and versatility of the SM-SPSA algorithm by solving a challenging problem based on web page ranking that Gurobi cannot solve due to the nonlinearity of the objective function. Unless otherwise stated, we use the centered mass heuristic (see Section~\ref{sec:infliction_points_remedies}).

\subsection{SM-SPSA versus Gurobi} \label{sec:numerical_experiments_SM-SPSA_vs_Gurobi}

In this section, we will compare the solution quality of the proposed SM-SPSA algorithm with the solution quality of solving the mathematical formulation of Section~\ref{sec:solver_approach} with Gurobi (version 10). If Gurobi has converged, then we know that this is also close to the theoretical optimal solution (up to a precision), making it an ideal benchmark for solution quality. As discussed in Section~\ref{sec:solver_approach}, Gurobi's time complexity to solve the mathematical formulation of Section~\ref{sec:solver_approach} is exponential. In contrast, the SM-SPSA algorithm has complexity $O(In^3)$. Therefore, it is expected that SM-SPSA scales better than Gurobi.

We generate 75 random network instances of 5 to 50 nodes according to the sampling method of Appendix~\ref{Appendix_sampling_networks}. The randomly generated binary adjustment matrices have zeroes on their diagonals to avoid trivial optimizations. The objective is to maximize the stationary distribution of a randomly chosen single node within each network. We run Gurobi until it has converged and run SM-SPSA for $750\times N^2$ iterations, as, in general, larger instances benefit from more iterations. The parameters used for both algorithms can be found in Appendix~\ref{Appendix_SMSPSA_Gurobi}.

The mean percentage difference of the objective values of the SM-SPSA algorithm compared to those of Gurobi is 1.77\%, which means that the SM-SPSA algorithm's solution quality is close to that of Gurobi. This is expected as the algorithm in \eqref{Standard_SA}, with $G $ taken as the expected value of the SM-SPSA update at $ \vartheta^{ (i)} $, is known to converge for $ \epsilon $ sufficiently small as $ i $ tends to $ \infty $ to an approximate stationary point of the objective; see Chapter~7 in \cite{KY}. This means that SM-SPSA will eventually find the (local) optimum. Appendix~\ref{Appendix_SMSPSA_Gurobi} presents a graph of the objective values found by SM-SPSA versus those found by Gurobi for all instances. 

Gurobi has an exponential time complexity. This indicates that the running time of a large network will blow up. In contrast, SM-SPSA also provides results for large networks in a reasonable time. To illustrate this, we generate a single instance with 500 nodes and run both algorithms for 24 hours. After these 24 hours, Gurobi does not have a solution, whereas the SM-SPSA algorithm achieves a significant increase in the objective value of $3486\%$ from 0.0020 to 0.0729. However, it should be noted that the SM-SPSA algorithm requires more iterations to converge to the optimal solution. Appendix~\ref{Appendix_SMSPSA_Gurobi} contains the used parameters for this experiment. The SM-SPSA result was obtained by not using the centered mass heuristic since larger networks have longer warm-up periods. In 24 hours, roughly 2.100.000 iterations could be performed, which is most likely too few to pass the warm-up period for such a large network. To verify this, we also ran the experiment using the centered mass heuristic, and then, indeed, a lower final objective value was obtained.







\subsection{Web Page Rank Optimization} \label{sec:numerical_experiments_hyperlink_optimization}

In this section, we will consider a real-life inspired web page ranking optimization problem that uses a challenging nonlinear objective function, which illustrates the versatility and applicability of the SM-SPSA algorithm.


The Internet is a large network where nodes represent web pages and edges represent hyperlinks between these web pages. When a user searches on a search engine, the search engine ranks the web pages according to their relevance. Algorithms such as Google PageRank \citep{Brin1998} are used to determine this web page ranking. These algorithms are often based on the number of incoming and/or outgoing links on a web page. Alternatively, the page ranking can also be based on the number of user clicks between web pages. In the latter case, more clicks lead to a higher page rank. A Markov chain can be derived by assigning weights proportional to the number of user clicks. The page rank of a web page then corresponds to the order statistic (read, ranking) of the stationary distribution of the Markov chain.

\begin{figure}[H]
    \centering
    \includegraphics[scale=0.4]{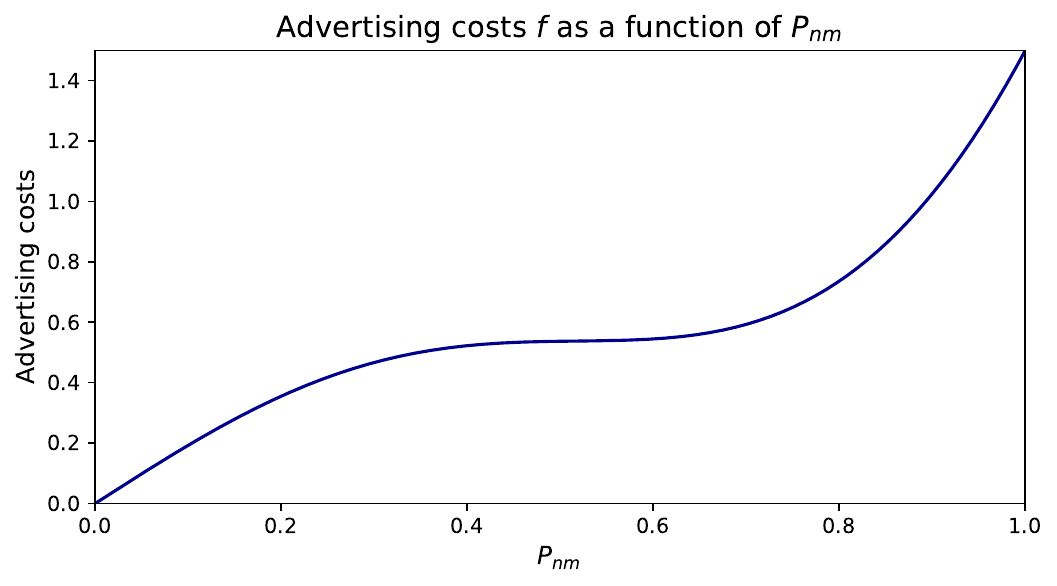}
    \caption{A graph of the function $f(P_{nm})$.}
    \label{fig:hyperlink_objective_function}
\end{figure}

Suppose that a web page owner wants to increase her page ranking by advertising on other web pages. As advertising requires investment, she wants to advertise on web pages that lead to the highest page ranking while minimizing costs. Let $m$ be the node of the web page owner, and let $J$ be the set of network nodes except $m$ ($J = S \setminus \{ m\}$). We consider the following objective function:

\begin{equation}\label{hyperlink_objective_1}
   \max \left (  \pi_m - \sum_{n\in J} \pi_n f(P_{nm})  \right ) ,
\end{equation}
where $f(P_{nm})$ represents the cost function dependent on the $P_{nm}$ entries, which can reflect, for example, that advertising costs for increasing web traffic from web page $n$ to $m$, as measured by $P_{nm}$. We multiply the costs $f$ by the stationary distribution of node $n$, $\pi_n$, to reflect that more popular websites charge higher advertising costs compared to unpopular websites. We consider the following nonlinear function for $f(P_{nm})$:

\begin{equation}\label{hyperlink_objective_2}
    f(P_{nm}) = 0.42\sin(1.5\pi P_{nm})+1.92(P_{nm})^3 .
\end{equation}
Figure~\ref{fig:hyperlink_objective_function} shows a graph of this function. The considered function reflects two phenomena common in advertising:
\begin{enumerate}
    \item Advertising costs are increasing in the number of advertisements.
    \item After buying a certain number of advertisements, there is a period of ``discount" where the costs of buying more advertisements remain constant. This is reflected by the period of near-stationarity of the costs between $\pm 0.4$ and $\pm 0.6$.
\end{enumerate}

\begin{figure}[H]
    \centering
    \includegraphics[scale=0.55]{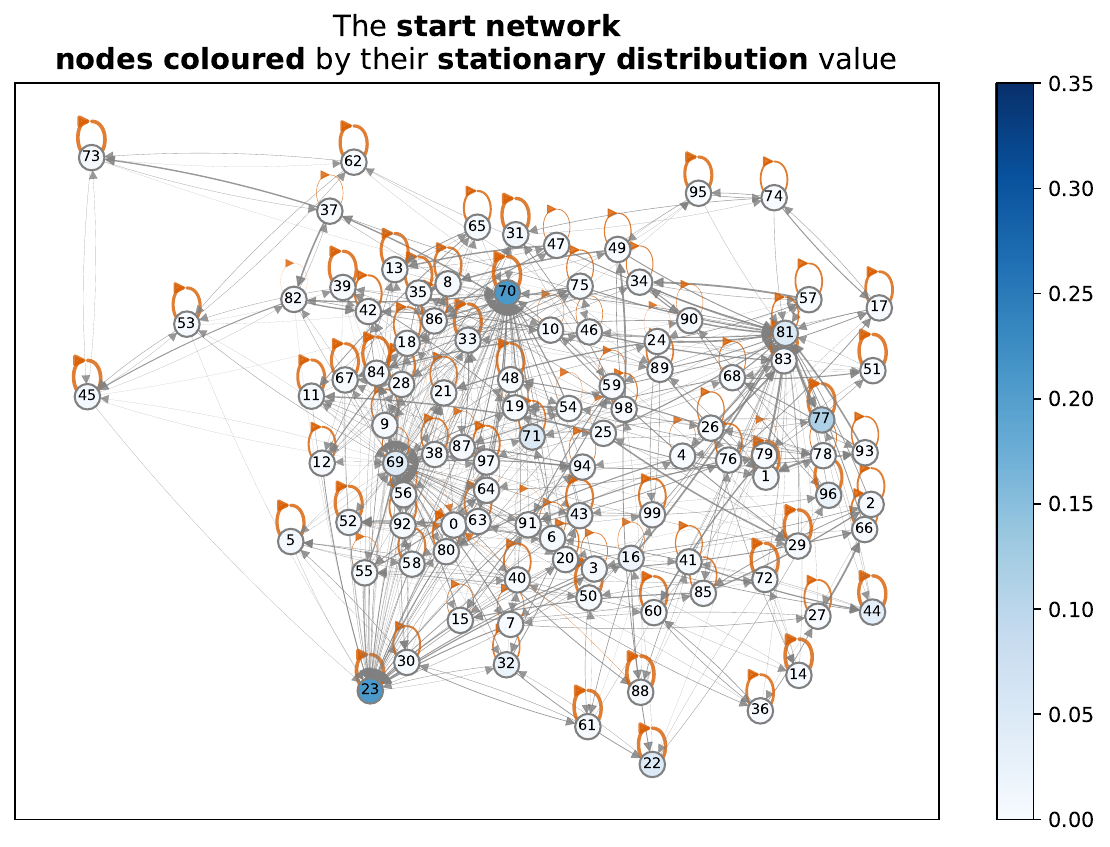}
    \caption{Strongly connected subcomponent of the Stanford web graph. Orange edges can be adjusted (non-existing adjustable edges are not visible). The darker the node color, the higher the value of the stationary distribution of that node.}
    \label{fig:hyperlink_network_before_optimization}
\end{figure}

We consider a strongly connected subcomponent (SCC) of 100 nodes of the Stanford web graph of 2002 \citep{StanfordWebgraph}. The web page owner can decide to buy advertisements on each web page $n\in J$, which means that the $C$ matrix has a value of 1 on all edges to node $m$. Furthermore, we assume that each web page $n \in J$ has some unused potential for advertisements, reflected in a self-loop for each node $n\in J$ that can also be adjusted. For each node, we randomly generate the mass of the self-loop from the uniform distribution from $0$ to $0.9$. The already existing edges in the strongly connected subcomponent also get a randomly generated mass assigned in such a way that the row-sums are equal to $1$ and thus a Markov chain is derived. We choose $m=0$, which means that node $0$ is the web page whose stationary distribution will be optimized. Figure~\ref{fig:hyperlink_network_before_optimization} shows the initial web-page network.

We run SM-SPSA with parameters as provided in Table~\ref{tab:hyperlink_SM-SPSA_parameters} in Appendix~\ref{Appendix_hyperlink}. Figure~\ref{fig:hyperlink_network_after_optimization} shows the network after optimization. Clearly, optimization resulted in a higher value of the stationary distribution for node $ 0 $ ($0.0010$ versus $0.3239$, respectively). 

A simple heuristic to evaluate the found solution is to run an optimization where only the stationary distribution of node $0$ is maximized without considering the costs. Table~\ref{tab:hyperlink_with_vs_without_costs_comparison} presents the stationary distribution of node $0$, the costs, and the final objective value according to \eqref{hyperlink_objective_1} and \eqref{hyperlink_objective_2} for the cases where the optimization is performed with and without taking the costs into account. As expected, the stationary distribution value of node $0$ is higher when the costs are not taken into account. However, the costs are considerable, which means that taking the costs into account results in a better final objective value. To conclude, this experiment illustrated that the SM-SPSA algorithm can be used for large networks with general nonlinear objective functions.

\begin{figure}[H]
    \centering
    \includegraphics[scale=0.55]{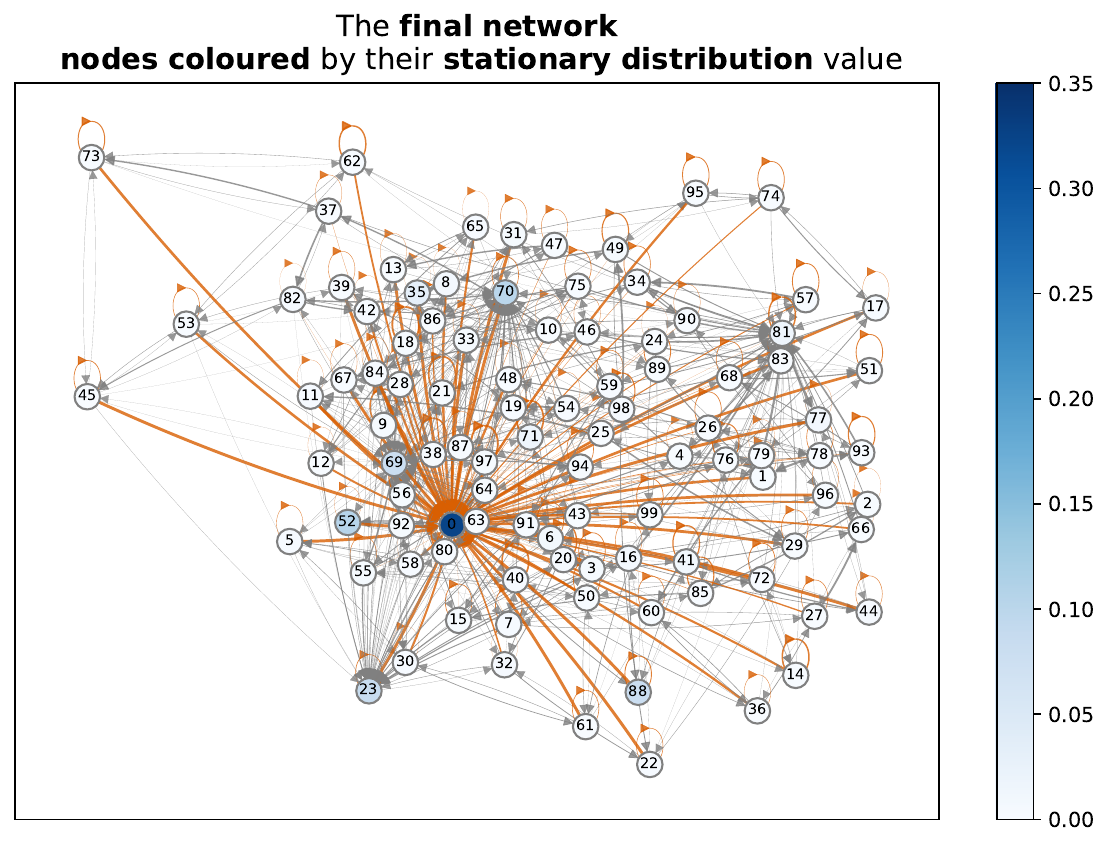}
    \caption{Network after maximizing the stationary distribution of node 0 while taking advertisement costs into account. Figure~\ref{fig:hyperlink_network_before_optimization} shows the original network.}
    \label{fig:hyperlink_network_after_optimization}
\end{figure}

\begin{table}[H]
    \centering
    \begin{tabular}{l|l|l}
         & Including costs & Excluding costs\\
        \hline
        Stationary distribution & 0.3239 & \cellcolor{green!10}0.3964\\
        Costs &\cellcolor{green!10}0.2798 & 0.4151\\
        Objective value &\cellcolor{green!10}0.0441 & -0.0187\\
    \end{tabular}
    \caption{A comparison between an optimization that does or does not take the costs into account.}
    \label{tab:hyperlink_with_vs_without_costs_comparison}
\end{table}

\section{Summary and Outlook on Further Research}\label{sec:discussion_conclusion}

In this paper, we introduced a framework for performing model-free optimization for objective functions over the stationary distribution of discrete-time Markov chains. The framework is built around the idea of combining first-order pseudo-gradient methods with a transformation to handle the stochasticity constraints of the problem. We used the framework to develop the SM-SPSA algorithm, which uses the well-known SPSA algorithm as pseudo-gradient method.

We identified and explained the non-standard behavior of the algorithm caused by infliction points, which arise due to small starting transition probabilities $(P_0)_{mn}$ in combination with our transformation. The existence of infliction points is independent of the applied (pseudo-)gradient method, which made it important to reduce this behavior. We a heuristic based on row-wise redistributing the values of the initial transition probabilities, called the centered mass heuristic, and provided statistical tests showing the positive influence of the heuristic on the convergence properties of the optimization algorithm.

Existing solvers can solve simple instances of the considered optimization problem, and we compared the performance of the SM-SPSA algorithm with the performance of the Gurobi solver for a wide variety of simple instances. These experiments showed that SM-SPSA  performs well (an average optimality gap of 1.77\%). We then applied the SM-SPSA algorithm to a challenging nonlinear problem to show the applicability and versatility of the developed algorithm.

To conclude, our SM-SPSA optimization algorithm provides an efficient and versatile method for optimizing functions over the stationary distribution of Markov chains. 
Our research also contributes to the field of stochastic approximation by identifying the phenomenon of ``infliction points" that arise when optimizing Markov chains and proposing a heuristic to reduce the impact of these infliction points.

The SM-SPSA algorithm can be applied to  large networks. However, the actual running time is greatly influenced by the implementation of the algorithm. In this paper, we used a standard implementation in Python, which was sufficient for running the experiments. However, there are many options to speed up the implementation of the SM-SPSA algorithm, such as using approximations for the stationary distribution, using one-measurement variants of SPSA \citep{Bhatnagar2013}, using parallel computing, or faster programming languages. Further research is on improving the running time of the SM-SPSA algorithm to solve (even) larger instances. Moreover, it is interesting to test the applicability of SM-SPSA on the optimization of other Markov chain measures beyond functions of the stationary distribution.


\bibliography{bibliography}

\appendix
\section{Appendix: Heuristic Statistical Tests}\label{Appendix_test}

Table~\ref{tab:infliction_point_parameters_no_heuristic} presents the used parameters for running SM-SPSA on the network in \eqref{eq:example_infliction_points}.

\begin{table}[H]
    \centering
    \begin{tabular}{l|l}
    Parameter     & Value  \\
    \hline
    Fixed gain size $\epsilon$     & 0.1 \\
    Number of iterations $I$ & 500.000 \\
    \end{tabular}
    \caption{The used parameters for running SM-SPSA on the network in \eqref{eq:example_infliction_points}.}
    \label{tab:infliction_point_parameters_no_heuristic}
\end{table}

\begin{figure}[H]
    \centering
    \centerline{
    \includegraphics[scale=0.45]{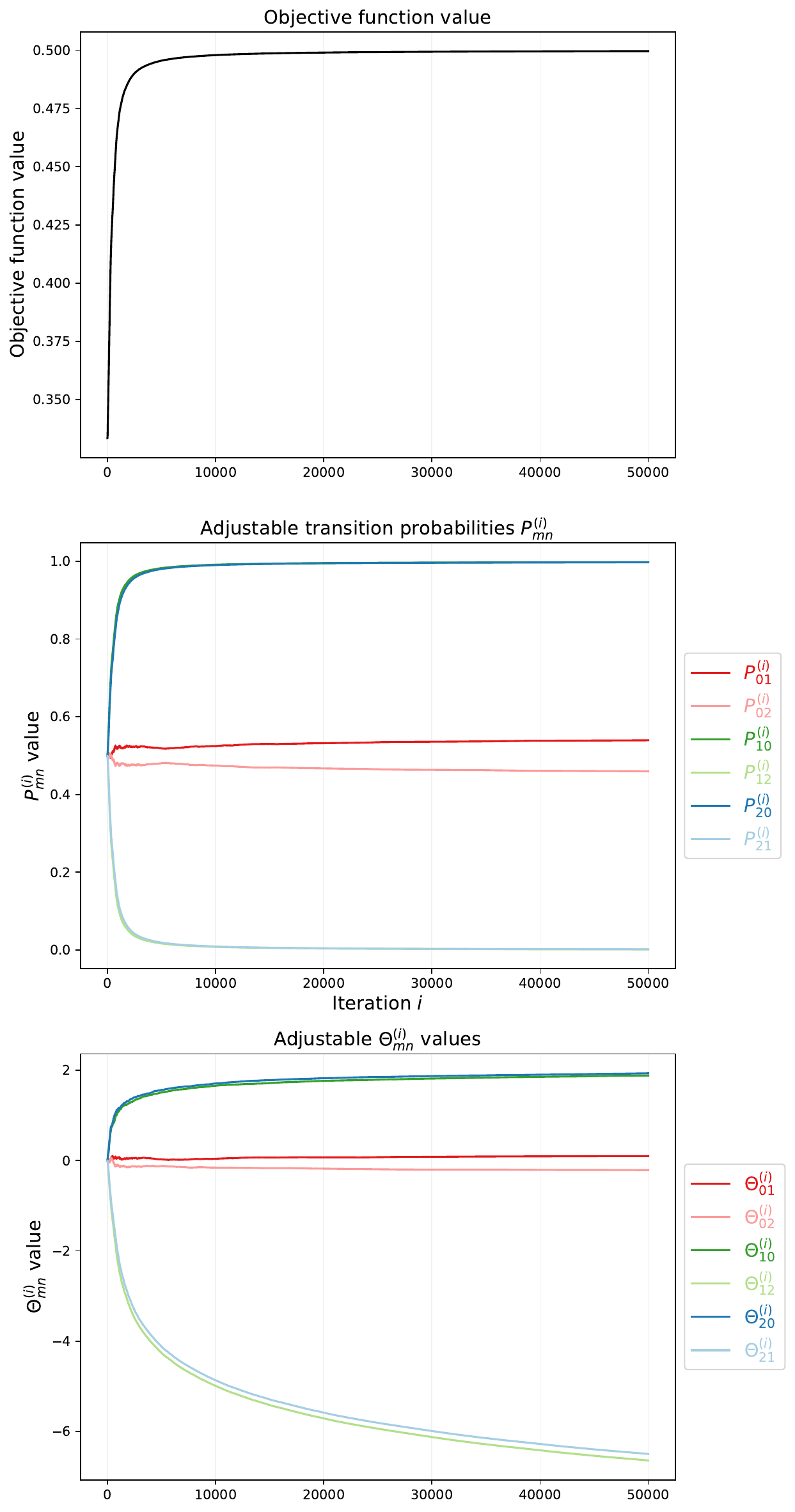}
    }
     \caption{Results of an optimization run for the network in \eqref{eq:example_infliction_points} with the centered mass heuristic.}
    \label{fig:infliction_points_example_centred_mass}
\end{figure}

Figure~\ref{fig:infliction_points_example_centred_mass} shows the three graphs similar to those in Figure~\ref{fig:infliction_points_example}, but then for an optimization with SM-SPSA and the centered mass heuristic. Table~\ref{tab:infliction_point_parameters_centred_mass} shows the used parameters. No infliction points can be observed, and the optimization converges to a similar objective function to that in Figure~\ref{fig:infliction_points_example} in only 50.000 iterations instead of 500.000 iterations; a significant improvement.

\begin{table}[H]
    \centering
    \begin{tabular}{l|l}
    Parameter     & Value  \\
    \hline
    Fixed gain size $\epsilon$     & 0.1 \\
    Number of iterations $I$ & 50.000 \\
    \end{tabular}
    \caption{The used parameters for running SM-SPSA on the network in \eqref{eq:example_infliction_points} with the centered mass heuristic.}
    \label{tab:infliction_point_parameters_centred_mass}
\end{table}

\begin{table}[H]
    \centering
    \begin{tabular}{l|l}
    Parameter  & Value  \\
    \hline
    Fixed gain size $\epsilon$     & 0.1 \\
    Number of iterations $I$ (medium, large) & 750.000 \\
    Number of iterations $I$  (small) & 100.000 \\
    \end{tabular}
    \caption{The used parameters for running SM-SPSA for the statistical tests.}
    \label{tab:statistical_experiment_parameters_SM-SPSA}
\end{table}

\section{Appendix: Sampling Random Networks}\label{Appendix_sampling_networks}

For the experiments, we need to use weighted networks whose transition matrices are stochastic. These networks are preferably generated randomly, so that experiments can be conducted with many different networks. Random networks can be sampled in multiple ways. In real-life settings, many networks contain weak links, i.e., edges with small weights. Therefore, we want to use a sampling method that generates networks with weak links. We use the following sampling method to sample a transition matrix of size $N\times N$:

\begin{enumerate}
    \item Generate an $N\times N$ matrix $Q$ where all entries are generated according to a uniform distribution on $[0,1]$.
    \item Generate a directed $N\times N$ Erdős-Rényi graph $P$, where an edge is present with probability $p_{ER}$.
    \item Combine the matrices $Q$ and $P$ according to $\alpha P+(1-\alpha)Q$, where $\alpha$ is a parameter between $0$ and $1$.
    \item Divide all entries of this new matrix row-wise by their row sum, such that a stochastic matrix is obtained.
\end{enumerate}

\noindent Unless otherwise stated, we use $\alpha=0.9$ and $p_{ER}=0.2$.

\section{Appendix: SM-SPSA versus Gurobi}\label{Appendix_SMSPSA_Gurobi}

Tables~\ref{tab:SM-SPSAvsGurobi_parameters_SM-SPSA} and \ref{tab:SM-SPSAvsGurobi_parameters_Gurobi} present the parameters used to run the experiments with the SM-SPSA and Gurobi algorithms, respectively. The built-in $FeasibilityTol$ parameter of Gurobi determines the precision/tolerance with which all constraints must be satisfied. The lower this value, the more precisely the constraints are met. However, the running time also increases when $FeasibilityTol$ decreases and numerical instability can also possibly become an issue. It should also be noted that the value of $\gamma$ should be greater than that of $FeasibilityTol$. In turn, $\gamma$ should be lower than the lowest transition probability that cannot be adjusted due to constraints \eqref{Mathematical_model_5} and \eqref{Mathematical_model_7} in Section~\ref{sec:solver_approach}. As the randomly generated networks can contain low transition probabilities and we use the Gurobi algorithm as a benchmark, we set $FeasibilityTol$ to the lowest possible value and take $\gamma$ one power of ten larger (see Table~\ref{tab:SM-SPSAvsGurobi_parameters_Gurobi}). The same parameters were used to run both algorithms for 24 hours, with the difference that the number of iterations of the SM-SPSA algorithm was set to a high value so that the algorithm could be automatically stopped after running for 24 hours.

\begin{table}[H]
    \centering
    \begin{tabular}{l|l}
    Parameter  & Value  \\
    \hline
    Fixed gain size $\epsilon$     & 0.1 \\
    Number of iterations $I$ & $750 \times N^2$ \\
    \end{tabular}
    \caption{The used parameters for running SM-SPSA for the SM-SPSA versus Gurobi experiment.}
    \label{tab:SM-SPSAvsGurobi_parameters_SM-SPSA}
\end{table}

\begin{table}[H]
    \centering
    \begin{tabular}{l|l}
    Parameter  & Value  \\
    \hline
    $FeasibilityTol$     & $1e-9$ \\
    $\gamma$ & $1e-8$\\
    \end{tabular}
    \caption{The used parameters for running Gurobi for the SM-SPSA versus Gurobi experiment.}
    \label{tab:SM-SPSAvsGurobi_parameters_Gurobi}
\end{table}

\begin{figure}[H]
    \centering
    \centerline{
    \includegraphics[scale=0.5]{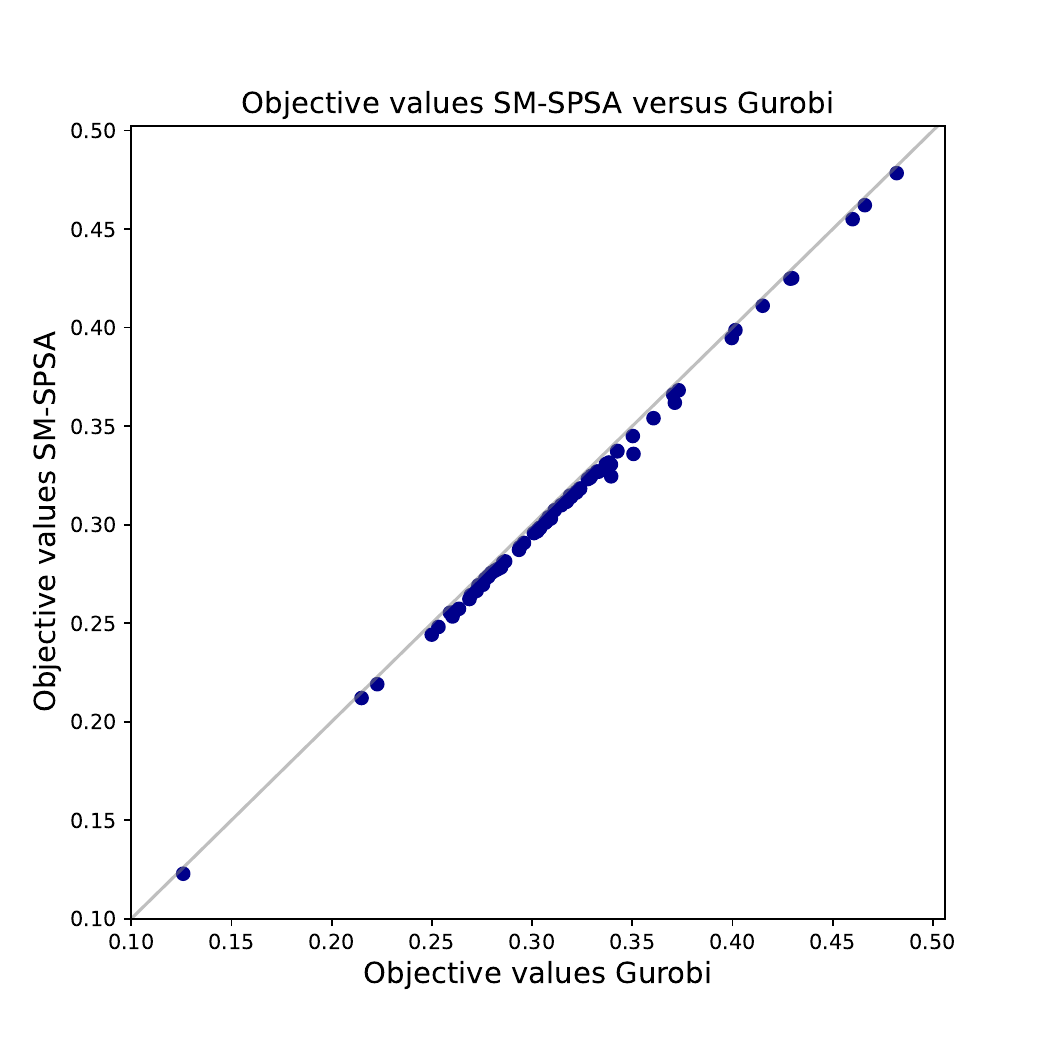}
    }
     \caption{The objective values of the SM-SPSA algorithm against those of Gurobi for the 75 generated instances.}
     \label{fig:SMSPSA_Gurobi_objective_values}
\end{figure}

Figure~\ref{fig:SMSPSA_Gurobi_objective_values} presents the objective values of the SM-SPSA algorithm versus those of Gurobi for the 75 generated instances. Clearly, the objective values of the SM-SPSA algorithm are close to those of Gurobi.

\section{Appendix: Web Page Rank Optimization}\label{Appendix_hyperlink}

\begin{table}[H]
    \centering
    \begin{tabular}{l|l}
    Parameter     & Value \\
    \hline
    Fixed gain size $\epsilon$     & 0.1 \\
    Number of iterations $I$ & 500.000\\
    \end{tabular}
    \caption{The used parameters for running SM-SPSA for the web page rank experiment.}
    \label{tab:hyperlink_SM-SPSA_parameters}
\end{table}

\newpage

\end{document}